\documentclass[preprint,12pt,a4paper,openright]{elsarticle}

\makeatletter
\def\ps@pprintTitle{%
	\let\@oddhead\@empty
	\let\@evenhead\@empty
	\def\@oddfoot{\reset@font\hfil\thepage\hfil}
	\let\@evenfoot\@oddfoot
}
\makeatother

\usepackage{mathptmx}	
\usepackage{enumitem}   
\usepackage[utf8]{inputenc}
\usepackage[english]{babel}
\usepackage{amsthm,amssymb,amsmath}
\usepackage{bigints}
\usepackage{multirow}
\usepackage{xcolor}
\usepackage{mathrsfs}
\usepackage{graphicx}
\usepackage{graphics} 
\usepackage{float} 
\usepackage{epsfig}
\usepackage{caption}
\usepackage{subcaption}
\usepackage[graphicx]{realboxes}
\usepackage{lscape} 

\usepackage{natbib}
\usepackage[colorlinks=true]{hyperref}
\hypersetup{
	colorlinks=true,
	linkcolor=blue,  
	citecolor=blue,
}

\newtheorem{cor}{Corollary}[section]

\newtheorem{thm}{Theorem}[section]

\title{On a Model for  Bivariate Left Censored Data}
\author{{G. Asha
		\footnote{Corresponding author email: asha@cusat.ac.in} 
		and Durga Vasudevan}\\
	{\textit{Department of Statistics}}\\
	{\textit{Cochin University of Science and Technology, Cochin-22}}\\
	{\textit{India.}}\\
	{Email: \textit{asha@cusat.ac.in, durgavvn10@gmail.com}}}

\begin{document}
	
	\begin{frontmatter}
		\begin{abstract}
			The lifetimes of subjects which are left-censored lie below a threshold value or a limit of detection. A popular tool used to handle left-censored data is the reversed hazard rate. In this work, we study the properties and develop characterizations of  a class of  distributions based on proportional reversed hazard rates used for analyzing left censored data. These characterizations are applied to simulate samples as well as analyze real data using distributions belonging to this class.
		\end{abstract}
	\end{frontmatter}
	
	\section{Introduction}
	A lifetime $Y$ associated with a subject is said to be left-censored if it is less than a censoring time $C_l$, which means that the event of interest has already occurred before that individual is considered in the study at time $C_l$. The exact lifetime $Y$ will be known if and only if $Y \geq C_l$. A left-censored data is represented by a pair of random variables $(Y,\delta)$ where $\delta=1$ if the event is observed and $0$ otherwise. Left-censored data have immense application in survival/reliability studies. They occur in life-test applications when a unit has failed at the time of its first inspection. They are also very common in bio-monitoring/environmental studies where observations could lie below a threshold value called limit of detection(LOD). Discarding the non-detected values for estimating parameters in a left-censored data set is a naive approach and many alternative techniques for handling left censored data have been proposed by many authors. 
	\par An easy way of approaching the left-censored datasets is substituting the non-detected values below the LOD with a constant such as the LOD, $\frac{LOD}{\sqrt(2)}$, or $\frac{LOD}{2}$. But when there is heavy censoring, these methods may produce error. The $\beta$ substitution method   computes a $\beta$-factor from the uncensored data in the dataset for adjusting the LOD that result in a near zero bias and lower root mean square error.  An algorithm for calculating the $\beta$-factor is provided in \citet{ganser2010accurate}.
	\par \citet{krishnamoorthy2009model} developed a model based imputation approach where observations below the LOD are randomly generated using the detected measurements. The parameters in the model are studied based on the detected observations and these imputed values. The data can then be analysed the using complete sample techniques along with adjustments to account for the imputation. 
	\par The conventional way is to estimate the parameters using the maximum likelihood estimation technique where the observed data contributes to the likelihood function through the probability density function (pdf), $f(y_i;\theta)$ and the left-censored data contributes through the cumulative distribution function (CDF), $F(y_i;\theta)$ as
	\begin{equation*}
		L \propto \prod_{i \in D} f(y_{i};\theta) \prod_{i \in C} F(y_{i};\theta),
	\end{equation*}
	where $D$ is the set of event times and $C$ is the set of left-censored observations. The non-parametric version of maximum likelihood estimation is the reverse Kaplan-Meier method.
	  The algorithm of this method is provided in (\citet{ware1976reanalysis}).

\par Another popular tool used to analyze the left-censored data sets is reversed hazard rate which was proposed by \citet{barlow1963properties} as a dual to the hazard rate. The concept of reversed hazard rate has been a subject of extensive study since it reappeared in \citet{keilson1982uniform} and \citet{block1998reversed}. \citet{ware1976reanalysis} reported the use of reversed hazard rate in the estimation of distribution function in the presence of left-censored observations.	Let $Y$ be an absolutely continuous random variable, $a=\text{inf}\{y~|F(y)>0\}$ and $b=\text{sup}\{y~|F(y)<1\}$. Then $I_{[a,b]}$ where $I_{[a,b]}=-\infty \leq a <b \leq \infty$ is the interval of support of $Y$ with the distribution function $F(y)$. The reversed hazard rate of $Y,$ denoted as $r(y),$ is defined for $y>a$ as,
	\begin{equation*}
		r(y)=\frac{f(y)}{F(y)},
	\end{equation*}
	where $f(y)$ and $F(y)$ are the probability density function and the distribution function of $Y$, respectively.
	\par  \citet{gupta1998modeling} proposed the proportional reversed hazards (PRH) model  which is expressed as,
	\begin{equation}\label{prhr}
		r(y)=\theta r_0(y),
	\end{equation}
	where $\theta>0$ and $r_0$ is the baseline reversed hazard rate. The corresponding distribution function is
	\begin{equation*}
		F(y)=[F_0(y)]^\theta,
	\end{equation*}
	where $F_0(y)$ is the baseline distribution function corresponding to $r_{0}(y).$ The model in \eqref{prhr} is helpful in the analysis of left-censored or right-truncated data (\citet{lawless2011statistical}). The PRH model has some extremely interesting properties. The parameter `$\theta$' is crucial in maintaining the structural properties of the baseline distribution. It is used to manage the skewness of the distribution. This is sometimes called the distortion Function too. For literature based on  distortion functions and PRH models we refer to \citet{arias2016new, ruggeri2021new}.  \citet{di2000some, mudholkar1993exponentiated, mudholkar1995exponentiated, gupta2007proportional, kizilaslan2017bayesian, popovic2021generalized}. When univariate ideas are extended to the bivariate setup, there could be  more than one extension since we need to capture the inherent dependence between the random variables. Hence there is more than one extension of the PRH model in the higher dimensions.  \citet{roy2002characterization} proposed a bivariate distribution function whose bivariate reversed hazard rates are locally proportional to the corresponding univariate reversed hazard rates and is given by,
	\begin{equation*}
		F(y_1,y_2)=F_{1}(y_1)F_{2}(y_2)exp\{-\gamma (log F_{1}(y_1))(log F_{2}(y_2))\},
	\end{equation*}
	for some $\gamma$. \citet{sankaran2008proportional} proposed a class of bivariate distributions given by,
	\footnotesize\begin{equation*}
		F(y_1,y_2)=(F_{1}(y_1))^{\alpha_{1}}(F_{2}(y_2))^{\alpha_{2}}exp\Big\{ \int_{y_1}^{\infty}\int_{y_2}^{\infty}(\beta_{1}m(u,v)-\beta_{2}\lambda_{1}(u|Y_{2} \leq v)\lambda_{2}(v|Y_{1} \leq u))dudv\Big\};~~\alpha_{i},\beta_{i}>0,~i=1,2,
	\end{equation*}
	\normalsize where $m(u,v)$ and $\lambda_{i}(y_{i}|Y_{3-i} \leq y_{3-i});~i=1,2$ are the bivariate reversed hazard rates as in \citet{bismi2005bivarite}, \citet{roy2002characterization} respectively. 
	\citet{kundu2010class} proposed a bivariate proportional reversed hazards model specified by
	\begin{equation}\label{kdf}
		F(y_1,y_2)=(F_{0}(y_1))^{\alpha_1}(F_{0}(y_2))^{\alpha_2}(F_{0}(z))^{\alpha_3},
	\end{equation}
	where $z=\text{min}\{y_1,y_2\}$. We refer this model as Bivariate Proportional Reversed Hazards Model 1 and denote as $ BPRHM1(F_{0},\alpha_1,\alpha_2,\alpha_{3})$. \citet{kundu2014multivariate} extended this model to the multivariate case in a similar manner. The models proposed by \citet{sankaran2008proportional}, \citet{kundu2010class} and \citet{kundu2014multivariate} have their marginals following univariate proportional reversed hazards model. \citet{durga2021proportional} extended the notion of PRH model to capture the inherent dependence where the status of the component affects the reversed hazard rate of the other component. They defined the bivariate density function of $(Y_1,Y_2)$ as,
	\footnotesize \begin{equation}\label{bdf}
		f_{Y_{1}, Y_{2}}\left(y_{1}, y_{2}\right)=\left\{\begin{array}{ll}{\theta_{1}^{\prime} \theta_{2} f_{0}\left(y_{1}\right) f_{0}\left(y_{2}\right)\left[F_{0}\left(y_{1}\right)\right]^{\theta_{1}^{\prime}-1}\left[F_{0}\left(y_{2}\right)\right]^{\theta_{1}+\theta_{2}-\theta_{1}^{\prime}-1}} & {\text { if } a<y_{1}<y_{2}<b} \\ {\theta_{1} \theta_{2}^{\prime} f_{0}\left(y_{1}\right) f_{0}\left(y_{2}\right)\left[F_{0}\left(y_{1}\right)\right]^{\theta_{1}+\theta_{2}-\theta_{2}^{\prime}-1}{\left[F_{0}\left(y_{2}\right)\right]^{\theta_{2}^{\prime}-1}}} & {\text { if } a<y_{2}<y_{1}<b}\end{array}\right.,
	\end{equation}
	\normalsize for $\theta_{i},\theta_i^{\prime}>0;~\text{i=1,2}$. We refer the model (\ref{bdf}) as Bivariate Proportional Reversed Hazards Model 2 and denote as $BPRHM2(F_{0},\theta_1,\theta_2,\theta_{1}^{\prime},\theta_{2}^{\prime})$. 
	The models in  (\ref{kdf}) and (\ref{bdf}) have a proportional reversed hazards model for $Y=\text{Max}(Y_1, Y_2)$. 
	\par In this paper we derive a class of distributions generalizing (\ref{kdf}) and (\ref{bdf}). This general class of distributions proposed includes many well studied bivariate distributions with interesting properties. 
For the general class of distributions proposed, we derive a characterization, appealing to which, we can develop a simple univariate procedure in lieu of complicated bivariate goodness of fit procedures. The outline of the paper is as follows.

	\par In Section \ref{chrresults}, we propose a class of distributions (BPRHM) based on a general functional equation. This class of distribution enjoys a PRH model for the distribution of component wise maxima. Characterizations of this  bivariate class of distributions, based on  functional equations, conditional mean and conditional variance are also derived in this section. In Section \ref{application}, we illustrate a simulation procedure for members of the BPRHM based on the characterizations developed in Section  \ref{chrresults}. Finally in Section \ref{realdata}, we conclude by illustrating the procedure for the American Football dataset of the National Football League [\citet{csorgo1989testing}] which was fitted using the \textit{BPRHM1} model (\ref{kdf}) with a Weibull baseline in \citet{kundu2010class}. Our findings are consistent with that in  \citet{kundu2010class}.


	
\section{Proposed Model and Characterization Properties} \label{chrresults}
	Consider a bivariate random pair $\textbf{Y}=(Y_1,Y_2)$ where $Y=\text{Max}\{Y_{1},Y_{2}\}$ has a univariate proportional reversed hazard (\textit{PRH}) rate with proportionality parameter $\theta$. Suppose $(Y_1,Y_2)$ has the cumulative distribution function $F(y_1,y_2)$  with support $I_{[a,b]}xI_{[a,b]}.$ Also assume that $F(y_1, y_2)$can be written in the form $F(y_1,y_2)=G(F_{0}(y_1),F_{0}(y_2))$, where $G$ and $F_0$, are such that,
	\begin{align}\label{cond}
		\begin{split}
		F(y_1,\infty)=&G(F_0(y_1),F_0(\infty))=G(F_0(y_1),1)=F_{Y_1}(y_1)\\
		F(\infty,y_2)=&G(F_0(\infty),F_0(y_2))=G(1,F_0(y_2))=F_{Y_2}(y_2)\\
		F(\infty,\infty)=&G(F_0(\infty),F_0(\infty))=G(1,1)=1.
		\end{split}
	\end{align} 
	Note that the function $G$ cannot be a copula since $G(1,F_0)=G(F_0,1) \ne F_0.$
	\begin{thm} \label{ch1}
		Under the conditions in \eqref{cond}, the distribution function $F(y_1,y_2)=G(F_0(y_1),F_0(y_2))$ satisfies the functional equation,
		\begin{equation} \label{gcrif}
			G(F_0(y)F_0(y_1),F_0(y)F_0(y_2))=G(F_0(y_1),F_0(y_2))G(F_0(y),F_0(y))
		\end{equation}
		for all $0<F_0(y),F_0(y_1),F_0(y_2)<1$, if and only if it is of the form,
		\begin{equation}\label{gcronly}
			G(F_0(y_1),F_0(y_2))=F\left(y_{1}, y_{2}\right)=\left\{\begin{array}{ll}{[F_0(y_1)]^{\theta} G \Big (1,\frac{F_0(y_2)}{F_0(y_1)} \Big )~~~\text{ if }~y_{1}\geq y_{2} }\\
				{[F_0(y_2)]^{\theta} G \Big (\frac{F_0(y_1)}{F_0(y_2)},1 \Big )~~~\text{ if }~y_{1}\leq y_{2} }\end{array}\right.
		\end{equation}
		for some y, $\theta>0$.
	\end{thm}
	
	\begin{proof}
From the conditions in \eqref{cond} and $y_1=y_2=y$
		\begin{equation*}
			G(F_0(y),F_0(y))=F(y,y)=[F_0(y)]^{\theta}
		\end{equation*}
		for some $\theta>0$. Therefore, \eqref{gcrif} can be written as,
		\begin{equation*}
			G(F_0(y)F_0(y_1),F_0(y)F_0(y_2))=G(F_0(y_1),F_0(y_2))[F_0(y)]^{\theta}.
		\end{equation*}
		Hence, for $y_1 \geq y_2$,
		\begin{align*}
			G(F_0(y_1),F_0(y_2))=&G \Big (F_0(y_1)F_0(\infty),F_0(y_1)\frac{F_0(y_2)}{F_0(y_1)} \Big) \\
			=&[F_0(y_1)]^{\theta} G \Big (1,\frac{F_0(y_2)}{F_0(y_1)} \Big).
		\end{align*}
		Similarly proceeding for $y_1 \leq y_2$, we arrive at \eqref{gcronly}.\\\\
		Conversely, 
		\normalsize from (\ref{gcronly}) for $y_1=y_2=y$, it follows that,
		\begin{equation*}
			G(F_0(y),F_0(y))=F(y,y)=[F_0(y)]^{\theta}.
		\end{equation*}
		\normalsize The result is now directly established by substituting  the distribution  \eqref{gcronly} in equation (\ref{gcrif}).
	\end{proof}

This representation of $F(y_1,y_2)$ as the product of the distribution of the maximum and marginals find applications that are mentioned in the next section. 
	
Assuming particular marginal distribution give rise to classes of distributions. 
	\begin{cor}\label{BPRHM}
		If the  marginals $F_{Y_1}(y_1)$ and $F_{Y_2}(y_2)$ have the  proportional reversed hazards model specified by,
		\begin{equation}\label{marginals1}
			F_{Y_i}(y_i)=[F_{0}(y_{i})]^{\alpha_{i}+\alpha_{3}},
		\end{equation}
		for $i=1,2$, then \eqref{gcronly} reduces to $BPRHM1(F_{0},\alpha_1,\alpha_2,\alpha_{3})$ (\citet{kundu2010class}) in \eqref{kdf} as,
		\begin{equation} \label{BPRHdf}
			F\left(y_{1}, y_{2}\right)=\left\{\begin{array}{ll}{(F_{0}(y_{1}))^{\alpha_{1}+\alpha_{2}+\alpha_{3}} \Big (\frac{F_{0}(y_{2})}{F_{0}(y_{1})} \Big )^{\alpha_{2}+\alpha_{3}} } & {\text { if } y_{1} \geq y_{2}} \\ { (F_{0}(y_{2}))^{\alpha_{1}+\alpha_{2}+\alpha_{3}} \Big (\frac{F_{0}(y_{1})}{F_{0}(y_{2})} \Big )^{\alpha_{1}+\alpha_{3}} } & {\text { if } y_{1} \leq y_{2}} \end{array}\right.
		\end{equation}
		for $\alpha_{i}>0;~\text{i=1,2,3}$.
	\end{cor}
	
	\begin{cor}\label{DPRHM}
		If the marginals $F_{Y_1}(y_1)$ and $F_{Y_2}(y_2)$ take the mixture form,
		\begin{equation}\label{marginals2}
			F_{Y_i}(y_i)=\frac{\theta_{3-i}}{\theta_i+\theta_{3-i}-\theta_i^{\prime}} [F_0(y_i)]^{\theta_i^{\prime}}+ \frac{\theta_i-\theta_i^{\prime}}{\theta_i+\theta_{3-i}-\theta_i^{\prime}}[F_0(y_i)]^{\theta_i+\theta_{3-i}}~~~\text{if}~a<y_{i}<b,
		\end{equation} 
		for $\theta_{1}+\theta_{2} \neq \theta_{i}^{\prime}$ and
		\begin{equation*}
			F_{Y_i}(y_i)=[F_{0}(y_{i})]^{\theta_{1}+\theta_{2}}[1-\theta_{3-i} ln (F_{0}(y_{i}))]~~~\text{if}~a<y_{i}<b,
		\end{equation*}
		for $\theta_{1}+\theta_{2}=\theta_{i}^{\prime}$, for $i=1,2$,
		then \eqref{gcronly} reduces to $BPRHM2(F_{0},\theta_1,\theta_2,\theta_{1}^{\prime},\theta_{2}^{\prime})$ specified by, (\citet{durga2021proportional}),
		\begin{equation}\label{DPRHdf}
			F\left(y_{1}, y_{2}\right)=\left\{\begin{array}{ll}{[F_0(y_1)]^{\theta_{1}+\theta_{2}} F_{Y_2} \Big [F_{0}^{-1} \Big( \frac{F_{0}(y_2)}{F_{0}(y_1)} \Big) \Big ]~~~\text{ if }~y_{1} \ge y_{2} }\\
				{[F_0(y_2)]^{\theta_{1}+\theta_{2}} F_{Y_1} \Big [F_{0}^{-1} \Big( \frac{F_{0}(y_1)}{F_{0}(y_2)} \Big) \Big ]~~~\text{ if }~y_{1} \le y_{2} }\end{array}\right.
		\end{equation}
		\normalsize for $\theta_{i},\theta_i^{\prime}>0;~\text{i=1,2}$.
	\end{cor}
Few interesting members of this class can be obtained by considering particular baseline distributions also.
\begin{cor}
	If the baseline distribution $F_{0}(y)=e^{c(y-b)};~c>0,~a<y \leq b,~b<\infty$, then \eqref{gcronly} becomes 
	\begin{equation}
		F\left(y_{1}, y_{2}\right)=\left\{\begin{array}{ll}{e^{c\theta (y_{1}-b)}F_{y_2}(b-(y_{1}-y_{2}))~~~\text{ if }~y_{1}\geq y_{2} }\\
			{e^{c\theta (y_{2}-b)}F_{y_1}(b-(y_{2}-y_{1}))~~~\text{ if }~y_{1}\leq y_{2} }\end{array}\right.
	\end{equation}
	for $c>0$, which is the general form of distribution satisfying the bivariate reversed lack of memory property (BRLMP),
	\begin{equation}\label{BRLMP}
		F(y_1,y_2)F(y,y)=F(0,0)F(y_1+y,y_2+y)
	\end{equation}
	for all $y_i$ and $y$ such that $a<y_i-y\leq y_i \leq b,~i=1,2$ (\citet{asha2007models}).
\end{cor}

	\begin{table}
		\caption{Functional equations satisfied by different bivariate distributions \label{tabfe}}
		\resizebox{1\textwidth}{!}{
			\begin{tabular}{llll}
				\hline
				& \textbf{Baseline distribution} & \hspace{3cm} & \textbf{Functional equation} \\ \hline
				&Reflected Weibull & & $F(y_1,y_2)F(y,y)=F(0,0)F(\sqrt{y_{1}^2+y^2},\sqrt{y_{2}^2+y^2})$ \\ & $F_{0}(y)=e^{-cy^2};~c>0,~y<0$ & & \\ \hline
				& Power Function & &  $F(y_1,y_2)F(y,y)=F(1,1)F \Big(\frac{yy_{1}}{b},\frac{yy_{2}}{b} \Big)$ \\ & $F_{0}(y)=\big(\frac{y}{b}\big)^c;~c>0,~0 \leq y < b,~b<\infty$& &          \\ \hline
				 & Inverse Exponential   &     &         $F(y_{1},y_{2})F(y,y)=F\Big (\frac{y_{1}y}{y_{1}+y},\frac{y_{2}y}{y_{2}+y} \Big )$ \\ & $F_0(y)=e^{-\frac{1}{y}};~~ y>0$ & &  \\ \hline
				& Exponential &   & $F(y_{1},y_{2})F(y,y)$ \\ & \multirow{2}{*}{ $F_0(y)=1-e^{-\lambda y};~~y>0, \lambda >0$} & &$=F\Big ( ln(e^{-\lambda y}+e^{-\lambda y_{1}} - e^{-\lambda (y_{1}+y)})^{-\frac{1}{\lambda}},$ \\ & & & $ln(e^{-\lambda y}+e^{-\lambda y_{2}} - e^{-\lambda (y_{2}+y)})^{-\frac{1}{\lambda}} \Big )$        \\ \hline
				& Inverse Weibull &   &  $F(y_{1},y_{2})F(y,y)=$ \\ & $	F_0(y)=e^{-y^{-\alpha}};~~y>0,~\alpha>0$ & & $F\Big ( (y_{1}^{-\alpha}+y^{-\alpha})^{-\frac{1}{\alpha}}, (y_{2}^{-\alpha}+y^{-\alpha})^{-\frac{1}{\alpha}} \Big )$             \\ \hline
				& Rayleigh &   &  $F(y_{1},y_{2})F(y,y)$ \\ & \multirow{2}{*}{ $F_0(y)=1-e^{-\lambda y^2};~~y>0$ }  & &$=F\Big ( \sqrt{ln(e^{-\lambda y^2}+e^{-\lambda y_{1}^2} - e^{-\lambda (y_{1}^2+y^2)})^{-\frac{1}{\lambda}}},$ \\ & & & $\sqrt{ln(e^{-\lambda y^2}+e^{-\lambda y_{2}^2} - e^{-\lambda (y_{2}^2+y^2)})^{-\frac{1}{\lambda}}} \Big )$       \\ \hline				
				& Weibull &   & $F(y_{1},y_{2})F(y,y)$ \\ & \multirow{2}{*}{$F_0(y)=1-e^{-\lambda y^{\beta}};~~y>0, \lambda >0, \beta>0$       } & &$=F\Big ( \big[ln(e^{-\lambda y^2}+e^{-\lambda y_{1}^2} - e^{-\lambda (y_{1}^2+y^2)})^{-\frac{1}{\lambda}}\big]^{\frac{1}{\beta}},$ \\ & & & $\big[ln(e^{-\lambda y^2}+e^{-\lambda y_{2}^2} - e^{-\lambda (y_{2}^2+y^2)})^{-\frac{1}{\lambda}}\big]^{\frac{1}{\beta}} \Big )$       \\ \hline   
			\end{tabular}
		}
	\end{table}

	More examples are provided in Table \ref{tabfe}. Hence the class \eqref{gcronly} is a rich class of distributions generally. Many class of distributions existing in literature can be obtained by virtue of the marginals and baseline distributions.
	\par Interesting characterization of class of distributions in \eqref{gcronly} based on the reversed hazard rate and moments of the distribution of $(Y_1,Y_2)$ are studied next. The concept of reversed hazard rate when extended to the higher dimensions takes into consideration various inherent dependence enjoyed between the components. This resulted in many extensions of this concept. See \citet{gurler1996bivariate, roy2002characterization}.

 \par \citet{roy2002characterization} defined the bivariate reversed hazard rate  of $\textbf{Y}= (Y_1, Y_2) $ as a two component vector,
 \begin{equation}\label{royrhr}
 	\lambda_{\textbf{Y}}(y_1,y_2)=\big (\lambda_{1}(y_1|Y_2 \leq y_2), \lambda_{2}(y_2|Y_1 \leq y_1) \big ),
 \end{equation}
 where $\lambda_{i}(y_i|Y_{3-i} \leq y_{3-i})= \lim_{\Delta y_i \to 0} \frac{P(y_i-\Delta y_i < Y_i \leq y_i|Y_1 \leq y_1, Y_2 \leq y_2)}{\Delta y_i}=\frac{\partial ln F(y_1,y_2)}{\partial y_i};~i=1,2$.
The definition in (\ref{royrhr}) uniquely determine the underlying distributions through

\normalsize 
\begin{equation}\label{roydf1}
	F(y_1,y_2)=exp \big[ -\int_{y_1}^{b_1} \lambda_{1}(u|Y_{2} \leq b_{2}) du - \int_{y_2}^{b_2}\lambda_{2}(v|Y_{1} \leq y_{1}) dv\big],
\end{equation}
or
\begin{equation}\label{roydf2}
	F(y_1,y_2)=exp \big[ -\int_{y_1}^{b_1} \lambda_{1}(u|Y_{2} \leq y_{2}) du - \int_{y_2}^{b_2}\lambda_{2}(v|Y_{1} \leq b_{1}) dv\big],
\end{equation}
respectively. Here, $\lambda_{1}(y_{1}|Y_{2} \leq b_{2})$ and $\lambda_{2}(y_2|Y_{1} \leq b_{1})$ are the marginal reversed hazard rates of $Y_{1}$ and $Y_{2}$, respectively.

 We propose a definition of  the bivariate reversed hazard rate of $\textbf{Y}=(Y_1,Y_2)$ as a dual of the conditional hazard rate of \citet{cox1972regression} as a three component vector,
\begin{equation}\label{bprhr}
	r_{\textbf{Y}}(y_1,y_2)=\big(r_{\textbf{Y}:10}(y)+r_{\textbf{Y}:20}(y),r_{\textbf{Y}:12}(y_1|y_2),r_{\textbf{Y}:21}(y_2|y_1)\big),
\end{equation}
where $r_{Y}(y)dy=(r_{\textbf{Y}:10}(y)+r_{\textbf{Y}:20}(y))dy$ with
\begin{equation*}
	r_{\textbf{Y}:i0}(y)= \lim_{\Delta y \to 0^{+}} \frac{P(y-\Delta y < Y_i \leq y|Y_1 \leq y, Y_2 \leq y)}{\Delta y}; \,y_i=y,~\text{i=1,2} 
\end{equation*}
and 
\begin{equation*}
	r_{\textbf{Y}:i3-i}(y_i|y_{3-i}) = \lim_{\Delta y_{i} \to 0^{+}} \frac{P(y_i-\Delta y_i < Y_i \leq y_i|Y_i \leq y_i, Y_{3-i} = y_{3-i})}{\Delta y_i};~y_i<y_{3-i},~\text{i=1,2}.
\end{equation*}
If $Y_{1}$ and $Y_{2}$ are independently distributed, then $r_{\textbf{Y}:i3-i}(y_i|y_{3-i})=r_{Y_{i}}(y_i)$ and $r_{\textbf{Y}:i0}(y)=r_{Y_{i}}(y);~i=1,2$. The bivariate reversed hazard rate $r_{\textbf{Y}}(y_1,y_2)$ in (\ref{bprhr}) uniquely determine the underlying distributions through the equation,
\footnotesize \begin{equation}
	f_{Y_{1}, Y_{2}}\left(y_{1}, y_{2}\right)=\left\{\begin{array}{ll}{r_{\textbf{Y}:20}(y_2) r_{\textbf{Y}:12}(y_1|y_2) exp \Big [- \int_{y_1}^{y_2} r_{\textbf{Y}:12}(u|y_2)du - \int_{y_2}^{b} \{r_{\textbf{Y}:10}(u)+r_{\textbf{Y}:20}(u)\}du \Big ]} & {\text { if } a<y_{1}<y_{2}<b} \\ {r_{\textbf{Y}:10}(y_1) r_{\textbf{Y}:21}(y_2|y_1) exp \Big [- \int_{y_2}^{y_1} r_{\textbf{Y}:21}(u|y_1)du - \int_{y_1}^{b} \{r_{\textbf{Y}:10}(u)+r_{\textbf{Y}:20}(u)\}du \Big ]} & {\text { if } a<y_{2}<y_{1}<b}\end{array}\right..
\end{equation}
\normalsize We propose few characterizations of the BPRHRM1 and BPRHRM2 based on the reversed hazard rates (\ref{royrhr}) and (\ref{bprhr}) using the following definitions given below.
\par Let  $\textbf{X}=(X_1, X_2)$, $Y=\text{Max}\{Y_1,Y_2\}$ and $X=\text{Max}\{X_1,X_2\}$. We define the following functions for all $y$, $y_{1}$ and $y_{2}$ such that $F_{X}(y)>0$ and $F_{X_{i}|X_{j}}(y_{i}|y_{j})>0;~i \neq j=1,2$ as,
	\begin{align*}
		\begin{split}
			A(y)&=-log F_{X}(y) \\
			A_{ij}(y_i|y_j)&=- ln F_{X_i|X_j=y_j}(y_i|y_j);~y_i<y_j;~i  \leq j=1,2 \\
			a_{Y;n}(y)&=E[A^{n}(Y)|Y_{1}<y,Y_{2}<y]=E[A^{n}(Y)|Y<y] \\
			a_{\textbf{Y};n}(y_{i}|y_{j})&=E[A_{ij}^{n}(Y_{i}|Y_{j})|Y_{i} \leq y_{i},Y_{j}=y_{j}];~y_i<y_j;~i \neq j=1,2,
		\end{split}
	\end{align*}
	where $n$ is a positive integer and $F_{X}(y)$ and $F_{\textbf{X}}(y_1,y_2)$ are the distribution functions of $X$ and $\textbf{X}$, respectively. Then $A^{n}(y)=[-log F_{X}(y)]^{n}$ and $\frac{d}{dy} [A^{n}(y)]=-nr_{X}(y)A^{n-1}(y)$. Similarly we get, $\frac{d}{dy_i}A_{ij}^{n}(y_i|y_j)=-n A_{ij}^{n-1}(y_i|y_j)r_{\textbf{X}:ij}(y_i|y_j);~i \neq j=1,2$. We further define the functions for all $y_{1}$ and $y_{2}$ such that $F_{\textbf{X}}(y_{1},y_{2})>0$,
	\begin{align*}
		\begin{split}
			B(y_1,y_2)&=- ln F_{\textbf{X}}(y_1,y_2) \\
			b_{Y_{i};n}(y_{1},y_{2})&=E_{Y_{i}}[B^{n}(Y_{1},Y_{2})|Y_{1} \leq y_{1},Y_{2} \leq y_{2}];~i=1,2,
		\end{split}
	\end{align*}
	where $n$ is a positive integer. Then we have, $\frac{d}{dy_i}B^{n}(y_1,y_2)=-n B^{n-1}(y_1,y_2)\lambda_{X_i}(y_i|X_2 \leq y_2);~i=1,2$. Based on the above definitions we propose few characterizations of the $BPRHM1(F_{0},\alpha_1,\alpha_2,\alpha_{3})$ and $BPRHM2(F_{0},\theta_1,\theta_2,\theta_{1}^{\prime},\theta_{2}^{\prime})$.
	\begin{thm} \label{DPRHmoments}
		The distribution function $F(y_1,y_2)$ of $(Y_1,Y_2)$  is distributed as BPRHM2 as in (\ref{DPRHdf}) if and only if,	
		\begin{enumerate}
			\item  for  marginals specified as in (\ref{marginals2}), $F(y_1,y_2)$ satisfies the functional equation  \eqref{gcrif}. \label{st22}
			\item  the bivariate reversed hazard rate $	r_{\textbf{Y}}(y_1,y_2)$ is of the form in \eqref{bprhr} given by,
			\begin{equation}\label{rhrform}
				r_{\textbf{Y}}(y_1,y_2)=	((\theta_1+\theta_2) r_0(y), \theta_{1}^{\prime}r_0(y_{1}), \theta_{2}^{\prime}r_0(y_{2}))
			\end{equation} 
			where $\theta_{i},\theta_i^{\prime}>0;~\text{i=1,2}$ and $r_{0}(y)$ is the reversed hazard rate corresponding to the baseline distribution.
%
			\item 
			\begin{enumerate}
				\item $a_{Y;n}(y)=A^{n}(y)+\Big( \frac{n}{\theta_{1}+\theta_{2}} \Big) a_{Y;n-1}(y)$ \\
				$V[A(Y)|Y<y]=\frac{1}{(\theta_{1}+\theta_{2})^{2}}$
				\item $a_{\textbf{Y};n}(y_{i}|y_{j}) = A_{ij}^{n}(y_{i}|y_{j})+\frac{n}{{\theta_{i}^{\prime}}} a_{\textbf{Y};n-1}(y_{i}|y_{j});~i,j=1,2,~i \neq j$ \\
				$V[A_{ij}(Y_{i}|Y_{j})|Y_{i} \leq y_i,Y_{j}=y_j]=\frac{1}{(\theta_{i}^{\prime})^{2}};~i,j=1,2,~i \neq j$,
				\newline for any real numbers $y,y_{i}$ and $y_{j}$  such that $y_{i}<y_{j}$ and $i,j=1,2;~i \neq j$.
			\end{enumerate} \label{st24}
		\end{enumerate}
	\end{thm}

	\begin{proof}
		\begin{enumerate}
			\item For marginals specified as in \eqref{marginals2}, the \textit{BPRHM2} is the only solution to the functional equation which follows from Corollary \ref{DPRHM}.
			\item If $	r_{\textbf{Y}}(y_1,y_2)$ is of the form \eqref{rhrform}, then the underlying bivariate density function of $(Y_1,Y_2)$ is given as in \eqref{bdf} which leads to the bivariate distribution function in \eqref{DPRHdf}. The converse is direct. 
			\item 
			\begin{enumerate}
				\item Observe that for $Y= \text{Max}\{Y_1,Y_2\}$, $F_{Y}(y)=[F_0(y)]^{\theta_1+\theta_2}$, a univariate  proportional hazards model with proportionality parameter $\theta_{1}+\theta_2.$ It now follows directly from  \citet{kundu2004characterizations} that 
				$E[A^{n}(Y)|Y<y]=A^{n}(y)+\Big( \frac{n}{\theta_{1}+\theta_{2}} \Big) E[A^{n-1}(Y)|Y<y]$ and $V[A(Y)|Y<y]=\frac{1}{(\theta_{1}+\theta_{2})^{2}}$.
				\item
		\par Assume that for $i,j=1,2;~i \neq j$ and $y_{i}<y_{j}$, $r_{\textbf{Y}:ij}(y_{i}|y_{j})=\theta_{i}^{\prime} r_{\textbf{X}:ij}(y_{i}|y_j);~\theta_{i}^{\prime}>0$.  Then,
		\begin{align*}
			a_{\textbf{Y};n}(y_{i}|y_{j})=&E[A_{ij}^{n}(Y_{i}|Y_{j})|Y_{i} \leq y_{i},Y_{j}=y_{j}] \\
			=&\frac{1}{F_{Y_i|Y_j=y_j}(y_i|y_j)} \int_{a}^{y_i} A_{ij}^{n}(t_i|y_j) f_{Y_i|Y_j=y_j}(t_i|y_j)dt_i \\
			=&\Big[\frac{A_{ij}^{n}(t_i|y_j)F_{Y_i|Y_j=y_j}(t_i|y_j)}{F_{Y_i|Y_j=y_j}(y_i|y_j)}\Big]_{a}^{y_i} + \int_{a}^{y_i} \frac{nA_{ij}^{n-1}(t_i|y_j)r_{\textbf{X}:ij}(t_i|y_j)F_{Y_i|Y_j=y_j}(t_i|y_j)}{F_{Y_i|Y_j=y_j}(y_i|y_j)}dt_i\\
			=&A_{ij}^{n}(y_i|y_j)+\frac{n}{\theta_{i}^{\prime} F_{Y_i|Y_j=y_j}(y_i|y_j)}\int_{a}^{y_i} A_{ij}^{n-1}(t_i|y_j)f_{Y_i|Y_j=y_j}(t_i|y_j)dt_i\\
			=&A_{ij}^{n}(y_i|y_j) + \Big( \frac{n}{\theta_{i}^{\prime}} \Big) a_{\textbf{Y};n-1}(y_i|y_j).
		\end{align*}
		Then, we have,
		\begin{align*}
			a_{\textbf{Y}:1}(y_i|y_j)=&A_{ij}(y_i|y_j)+\frac{1}{\theta_{i}^{\prime}},~\text{and}\\
			a_{\textbf{Y};2}(y_i|y_j)=&A_{ij}^{2}(y_i|y_j)+\frac{2}{\theta_{i}^{\prime}}a_{\underline{Y}:1}(y_i|y_j) \\
			=&A_{ij}^{2}(y_i|y_j)+\frac{2}{\theta_{i}^{\prime}} \Big (A_{ij}(y_i|y_j)+\frac{1}{\theta_{i}^{\prime}} \Big ) \\
			=&\Big [A_{ij}(y_i|y_j)+\frac{1}{\theta_{i}^{\prime}} \Big]^{2} + \frac{1}{(\theta_{i}^{\prime})^{2}}\\
			=& \big (a_{\textbf{Y}:1}(y_i|y_j) \big )^{2}+\frac{1}{(\theta_{i}^{\prime})^{2}}.
		\end{align*}
		Then,
		\begin{align*}
			V[A_{ij}(Y_i|Y_j)|Y_i \leq y_i, Y_j=y_j]=&E[A_{ij}^{2}(Y_i|Y_j)|Y_i \leq y_i, Y_j=y_j]-\big[ E[A_{ij}(Y_i|Y_j)|Y_i \leq y_i, Y_j=y_j] \big]^{2}\\
			=&a_{\textbf{Y};2}(y_i|y_j)-\big (a_{\textbf{Y}:1}(y_i|y_j) \big )^{2}\\
			=&\frac{1}{(\theta_{i}^{\prime})^{2}}.
		\end{align*}
		Conversely, suppose that
		\begin{eqnarray}
			a_{\textbf{Y};n}(y_i|y_j)=A_{ij}^{n}(y_i|y_j) + \Big( \frac{n}{\theta_{i}^{\prime}} \Big) a_{\textbf{Y};n-1}(y_i|y_j) \label{con1}\\ \text{and}~
			V[A_{ij}(Y_i|Y_j)|Y_i \leq y_i, Y_j=y_j]=\frac{1}{(\theta_{i}^{\prime})^{2}} \label{con2}
		\end{eqnarray}
		is true for $i,j=1,2;~i \neq j$. Then \eqref{con1} implies
		\footnotesize \begin{equation} \label{ch21}
			\int_{a}^{y_i} A_{ij}^{n}(t_i|y_j)f_{Y_i|Y_j=y_j}(t_i|y_j)dt_i = A_{ij}^{n}(y_i|y_j)F_{Y_i|Y_j=y_j}(y_i|y_j)+\Big( \frac{n}{\theta_{i}^{\prime}} \Big) \int_{a}^{y_i} A_{ij}^{n-1}(t_i|y_j)f_{Y_i|Y_j=y_j}(t_i|y_j)dt_i.
		\end{equation}
		\normalsize Differentiating both sides of \eqref{ch21} w.r.t $t_i$, we get,
		\begin{align*}
			&r_{\textbf{X}:ij}(y_i|y_j)A_{ij}^{n-1}(y_i|y_j)F_{Y_i|Y_j=y_j}(y_i|y_j)=\Big( \frac{1}{\theta_{i}^{\prime}} \Big) A_{ij}^{n-1}(y_i|y_j)f_{Y_i|Y_j=y_j}(y_i|y_j)\\
			&\implies A_{ij}^{n-1}(y_i|y_j) \Big[r_{\textbf{X}:ij}(y_i|y_j)F_{Y_i|Y_j=y_j}(y_i|y_j)-\Big( \frac{1}{\theta_{i}^{\prime}} \Big)f_{Y_i|Y_j=y_j}(y_i|y_j) \Big]=0.
		\end{align*}
		Since $A_{ij}^{n-1}(y_i|y_j) \neq 0$, we have, $r_{\textbf{Y}:ij}(y_i|y_j)=\theta_{i}^{\prime} r_{\textbf{X}:ij}(y_i|y_j)$.
		\par Now, \eqref{con2} implies
		\footnotesize \begin{align}
			\nonumber &\frac{1}{F_{Y_i|Y_j=y_j}(y_i|y_j)} \int_{a}^{y_i} A_{ij}^{2}(t_i|y_j)f_{Y_i|Y_j=y_j}(t_i|y_j)dt_i - \Big[\frac{1}{F_{Y_i|Y_j=y_j}(y_i|y_j)} \int_{a}^{y_i} A_{ij}(t_i|y_j)f_{Y_i|Y_j=y_j}(t_i|y_j)dt_i \Big]^{2}=\frac{1}{(\theta_{i}^{\prime})^{2}} \\
			\label{ch31} &\implies F_{Y_i|Y_j=y_j}(y_i|y_j) \int_{a}^{y_i} A_{ij}^{2}(t_i|y_j)f_{Y_i|Y_j=y_j}(t_i|y_j)dt_i - \Big[ \int_{a}^{y_i} A_{ij}(t_i|y_j)f_{Y_i|Y_j=y_j}(t_i|y_j)dt_i\Big]^{2}=\Big( \frac{F_{Y_i|Y_j=y_j}(y_i|y_j)}{\theta_{i}^{\prime}} \Big)^{2}.
		\end{align}
		\normalsize Differentiating both sides of \eqref{ch31} w.r.t $y_i$, we get,
		\footnotesize \begin{align} \label{ch32}
			\nonumber A_{ij}^{2}(y_i|y_j)F_{Y_i|Y_j=y_j}(y_i|y_j) + \int_{a}^{y_i} A_{ij}^{2}(t_i|y_j)f_{Y_i|Y_j=y_j}(t_i|y_j)dt_i -2A_{ij}(y_i|y_j)&\int_{a}^{y_i} A_{ij}(t_i|y_j)f_{Y_i|Y_j=y_j}(t_i|y_j)dt_i\\
			&=\frac{2 F_{Y_i|Y_j=y_j}(y_i|y_j)}{(\theta_{i}^{\prime})^{2}}.
		\end{align}
		\normalsize Differentiating both sides of \eqref{ch32} w.r.t $y_i$, we get,
		\footnotesize \begin{align}
			\nonumber A_{ij}(y_i|y_j)r_{\textbf{X}:ij}(y_i|y_j)\frac{d}{dy_i}F_{Y_i|Y_j=y_j}(y_i|y_j) - r_{\textbf{X}:ij}(y_i|y_j)\int_{a}^{y_i} A_{ij}(t_i|y_j)f_{Y_i|Y_j=y_j}(t_i|y_j)dt_i=&-\frac{f_{Y_i|Y_j=y_j}(y_i|y_j)}{(\theta_{i}^{\prime})^{2}}\\
			\label{ch33} A_{ij}(y_i|y_j)F_{Y_i|Y_j=y_j}(y_i|y_j)- \int_{a}^{y_i} A_{ij}(t_i|y_j)f_{Y_i|Y_j=y_j}(t_i|y_j)dt_i=&-\frac{f_{Y_i|Y_j=y_j}(y_i|y_j)}{(\theta_{i}^{\prime})^{2} r_{\textbf{X}:ij}(y_i|y_j)}.
		\end{align}
		\normalsize Differentiating both sides of \eqref{ch33} w.r.t $y_i$, we get,
		\begin{equation} \label{ch34} 
			(\theta_{i}^{\prime})^{2}F_{Y_i|Y_j=y_j}(y_i|y_j)=\frac{1}{r_{\textbf{X}:ij}(y_i|y_j)}\frac{d}{dy_i} \Big[ \frac{f_{Y_i|Y_j=y_j}(y_i|y_j)}{r_{\textbf{X}:ij}(y_i|y_j)}\Big].
		\end{equation}
		Now, \eqref{ch34} can be written in the form of a second order homogeneous differential equation as,
		\footnotesize \begin{equation}
			\frac{d^2}{dy_{i}^{2}} \Big[ F_{Y_i|Y_j=y_j}(y_i|y_j) \Big] - \frac{\frac{d}{dy_i}r_{\textbf{X}:ij}(y_i|y_j)}{r_{\textbf{X}:ij}(y_i|y_j)}\frac{d}{dy_i}\Big[ F_{Y_i|Y_j=y_j}(y_i|y_j) \Big] - (\theta_{i}^{\prime})^{2}r_{\textbf{X}:ij}^{2}(y_i|y_j)\Big[ F_{Y_i|Y_j=y_j}(y_i|y_j) \Big] =0
		\end{equation}
		\normalsize The general solution of this homogeneous problem is,
		\begin{align*}
			F_{Y_i|Y_j=y_j}(y_i|y_j)=&c_{1}e^{\theta_{i}^{\prime} A_{ij}(y_i|y_j)}+c_{2}e^{-\theta_{i}^{\prime} A_{ij}(y_i|y_j)}\\
			=&c_{1}[ F_{X_i|X_j=y_j}(y_i|y_j)]^{-\theta_{i}^{\prime}}+c_{2}[ F_{X_i|X_j=y_j}(y_i|y_j)]^{\theta_{i}^{\prime}}
		\end{align*}
		where $c_{1}$ and $c_{2}$ are two arbitrary constants. Now, as $y_{i} \rightarrow a$, we have $F_{Y_i|Y_j=y_j}(y_i|y_j) \rightarrow 0$ and $F_{X_i|X_j=y_j}(y_i|y_j) \rightarrow 0$ and $y_i \rightarrow b$, we have $F_{Y_i|Y_j=y_j}(y_i|y_j) \rightarrow 1$ and $F_{X_i|X_j=y_j}(y_i|y_j) \rightarrow 1$. Using these conditions we obtain $c_{1}=0$ and $c_{2}=1$. Hence, we have,
		\begin{equation*}
			F_{Y_i|Y_j=y_j}(y_i|y_j)=[F_{X_i|X_j=y_j}(y_i|y_j)]^{\theta_{i}^{\prime}} \implies r_{\textbf{Y}:ij}(y_i|y_j)=\theta_{i}^{\prime} r_{\textbf{X}:ij}(y_i|y_j).
		\end{equation*}
		\par If $X_{1}$ and $X_{2}$ are independently distributed, then $r_{\textbf{X}:ij}(y_i|y_j)=r_{X_{i}}(y_i)$ for $i,j=1,2;~i \neq j$. That is, for any real numbers $y_i$ and $y_j$ such that $y_i<y_j$ and $i,j=1,2;~i \neq j$, $r_{\textbf{Y}:ij}(y_i|y_j)=\theta_{i}^{\prime} r_{X_{i}}(y_i);~\theta_{i}^{\prime}>0$ if and only if either of the following is true,
		\begin{enumerate}
			\item $a_{\textbf{Y};n}(y_i|y_j)=A_{ij}^{n}(y_i|y_j) + \Big( \frac{n}{\theta_{i}^{\prime}} \Big) a_{\textbf{Y};n-1}(y_i|y_j),$ or \\ $E[A_{ij}^{n}(Y_{i}|Y_{j})|Y_{i} \leq y_i,Y_{j}=y_j] = A_{ij}^{n}(y_{i}|y_{j})+\frac{n}{{\theta_{i}^{\prime}}} E[A_{ij}^{n-1}(Y_{i}|Y_{j})|Y_{i} \leq y_i,Y_{j}=y_j]$
			\item $V[A_{ij}(Y_{i}|Y_{j})|Y_{i} \leq y_i,Y_{j}=y_j]=\frac{1}{\theta_{i}^{\prime 2}};~i,j=1,2,~i \neq j$,
		\end{enumerate}
	where $r_{\textbf{Y}:ij}(y_i|y_j)$ is the reversed hazard rate of $Y_{i}|Y_{j}$. Hence the proof.
		\end{enumerate}
	\end{enumerate}
\end{proof}

\begin{thm}\label{BPRHmoments}
	The distribution function $F(y_1,y_2)$ of $(Y_1,Y_2)$  is distributed as BPRHM1 as in (\ref{BPRHdf}) if and only if,
	\begin{enumerate}
		\item the distribution function $F(y_1,y_2)$ of $(Y_1,Y_2)$ with marginals specified in Corollary \ref{BPRHM} satisfy the functional equation in \eqref{gcrif}. \label{st31}
		\item $(Y_1,Y_2)$ has the bivariate reversed hazard rate as defined in \eqref{royrhr} given by,
		\begin{equation*}
			\lambda_{Y_{i}}(y_i|Y_{3-i} \leq y_{3-i})=\left\{\begin{array}{ll}{(\alpha_{i}+\alpha_{3})r_0(y_{i})~~~~~~~~~~~;~y_{i} < y_{3-i} }\\
			{\alpha_{i}r_0(y_{i})~~~~~~~~~~~~~~~~~~~~~~;~y_{i} > y_{3-i} }\\
			{(\alpha_{1}+\alpha_{2}+\alpha_{3})r_0(y)~~~;~y_{i}=y_{3-i}=y}\end{array}\right.,
		\end{equation*}
	for $i=1,2$ and $r_{0}(y)$ is the reversed hazard rate corresponding to the baseline distribution. \label{st32}
	\item \begin{enumerate}
		\item $a_{Y;n}(y)=A^{n}(y)+\Big( \frac{n}{\alpha_{1}+\alpha_{2}+\alpha_{3}} \Big) a_{Y;n-1}(y)$ \\
		$V[A(Y)|Y<y]=\frac{1}{(\alpha_{1}+\alpha_{2}+\alpha_{3})^{2}}$
		\item $b_{Y_{1};n}(y_{1},y_{2}) = B^{n}(y_{1},y_{2})+\frac{n}{{\alpha_{1}+(j-1)\alpha_{3}}} b_{Y_{1};n-1}(y_{1},y_{2})$ \\
		$b_{Y_{2};n}(y_{1},y_{2}) = B^{n}(y_{1},y_{2})+\frac{n}{{\alpha_{2}+(i-1)\alpha_{3}}} b_{Y_{2};n-1}(y_{1},y_{2})$ \\
		$V_{Y_1}[B(Y_{1},Y_{2})|Y_{1} \leq y_1,Y_{2} \leq y_2]=\frac{1}{(\alpha_{1}+(j-1)\alpha_{3})^2}$\\
		$V_{Y_2}[B(Y_{1},Y_{2})|Y_{1} \leq y_1,Y_{2} \leq y_2]=\frac{1}{(\alpha_{1}+(i-1)\alpha_{3})^2}$,
		\newline for any real numbers $y,y_{i}$ and $y_{j}$  such that $y_{i}<y_{j}$ and $i,j=1,2;~i \neq j$.
	\end{enumerate} \label{st33}
	\end{enumerate}
\end{thm}

\begin{proof}
	The equivalency of statements \ref{st31} and \ref{st32} can be proved using the unique representation of distribution function by \citet{roy2002characterization} given in \eqref{roydf1} or \eqref{roydf2}. Also,
	part (a) of statement \ref{st33} follows similarly as the part (a) of statement \ref{st24} in Theorem \ref{DPRHmoments}. Now, assume that for $i=1,2$,
	\begin{equation*}
		\lambda_{Y_{i}}(y_i|Y_{3-i} \leq y_{3-i})=\left\{\begin{array}{ll}{(\alpha_{i}+\alpha_{3})\lambda_{X_{i}}(y_i|X_{3-i} \leq y_{3-i})~~~~~~~~~~~;~y_{i} < y_{3-i} }\\
		{\alpha_{i}\lambda_{X_{i}}(y_i|X_{3-i} \leq y_{3-i})~~~~~~~~~~~~~~~~~~~~~~;~y_{i} > y_{3-i} }\end{array}\right..
	\end{equation*}
	Then, the results in part (b) of statement \ref{st33} follows proceeding similarly as in the proof of part (b) of statement \ref{st24} in Theorem \ref{DPRHmoments}.
\end{proof}

\section{Simulations }\label{application}
	\par The following algorithm based on theorem \ref{ch1} is used for generating the random variables from \textit{BPRHM} model.  The samples generated are then verified using the bivariate K-S test (\citeauthor*{justel1997multivariate}). The algorithm is as follows.
 	\subsection*{\normalfont Algorithm:}
 	\begin{itemize}
 		\item Generate  five independent uniform $(0,1)$ random variables $U_i,~i=1,2,3,4,5$.
 		\item Decide on censoring percentage $p$.
 		\item For a prefixed censoring percentage $p$ in the population, generate the censoring times $c_{1}=z_{1}U_{4}$ and $c_{2}=z_{2}U_{5}$, where $z_{1}$ and $z_{2}$ are derived by solving $P[0 \leq Y_{i} \leq c_{i}, 0 \leq c_{i} \leq z_{i}]=p;~i=1,2$.
 		\item  For a given baseline $F_0$, If $U_{1} \geq P(Y_{1} > Y_{2})$, $t_1=F_0^{-1}(U_{2}^{\frac{1}{\theta}})$ and $t_2=F_0^{-1}[F_{0}(t_1)F_0(F_{Y_2}^{-1}(U_{3}))]$.
 		\item If $U_{1} < P(Y_{1} > Y_{2})$, $t_2=F_0^{-1}(U_{2}^{\frac{1}{\theta}})$ and $t_1=F_0^{-1}[F_{0}(t_2)F_0(F_{Y_1}^{-1}(U_{3}))]$.
 		\item $y_{1}=\text{Max}\{t_{1},c_{1}\}$ and $y_{2}=\text{Max}\{t_{2},c_{2}\}$
 		\item Repeat n times for a sample of size $n$.
 	\end{itemize}
\par Samples of size $n=100$ from the $BPRHM1(W(1.5,1.2),1.3,1.2,1.0)$ model with Weibull $(W(\lambda,\beta))$ as baseline distribution whose CDF is $F(y;\lambda,\beta)=1-e^{-\lambda y^{\beta}};~y>0,\lambda>0,\beta>0$ were generated. The bivariate goodness of fit test  (\citet{justel1997multivariate}) is obtained as $p>0.10$ confirming the sample. Alternatively, univariate K-S tests for goodness of fit for $Y=\text{Max}\{Y_{1},Y_{2}\}$, the marginals, $Y_{1}$ and $Y_{2}$ gave consistent results in tune with Corollary \ref{BPRHM}. The results, are presented in Table \ref{tabks}. 
\par Bivariate samples from $BPRHM2(IW(1.2),1.2,1.4,1.6,1.8)$ model with Inverse Weibull $(IW(\beta))$ as baseline distribution whose CDF is given by, $F(y;\beta)=e^{-y^{-\beta}};~y>0,\beta>0$ was generated and illustrated similarly.
 Results under $20\%$ censoring is presented in Table \ref{tabks1}.

\begin{table}
	\centering
	\caption{Goodness of fit test for different baseline distributions (Uncensored samples) \label{tabks}}
	\resizebox{1\columnwidth}{!}{%
	\begin{tabular}{llccccc}
		\hline
		\multirow{2}{*}{Model}                                                                     & \multirow{2}{*}{\begin{tabular}[c]{@{}l@{}}Baseline\\ Distribution\end{tabular}}   & \multicolumn{1}{l}{\multirow{2}{*}{\begin{tabular}[c]{@{}l@{}}Bivariate\\ K-S Statistic\end{tabular}}} & \multicolumn{1}{l}{\multirow{2}{*}{P value}} & \multirow{2}{*}{Variable} & \multicolumn{1}{l}{\multirow{2}{*}{\begin{tabular}[c]{@{}l@{}}K-S\\ statistic\end{tabular}}} & \multicolumn{1}{l}{\multirow{2}{*}{P value}} \\ 
		&                                                                                    & \multicolumn{1}{l}{}                                                                                   & \multicolumn{1}{l}{}                         &                           & \multicolumn{1}{l}{}                                                                         & \multicolumn{1}{l}{}                         \\ \hline
		\multirow{3}{*}{\begin{tabular}[c]{@{}l@{}}\textit{BPRHM1}\\ (W(1.5,1.2),1.3,1.2,1.0)\end{tabular}} & \multirow{3}{*}{\begin{tabular}[c]{@{}l@{}}Weibull\\ \end{tabular}}   & \multirow{3}{*}{0.1557}                                                                                & \multirow{3}{*}{\textgreater 0.1000}                      & Max$\{Y_1,Y_2\}$          & 0.0484                                                                                       & 0.973                                        \\
		&                                                                                    &                                                                                                        &                                              & Marginal of $Y_1$         & 0.1096                                                                                       & 0.1812                                       \\
		&                                                                                    &                                                                                                        &                                              & Marginal of $Y_2$         & 0.1051                                                                                       & 0.2195                                       \\ \hline
		\multirow{3}{*}{\begin{tabular}[c]{@{}l@{}}\textit{BPRHM1}\\ (E(1.2),1.3,1.2,1.0)\end{tabular}}     & \multirow{3}{*}{\begin{tabular}[c]{@{}l@{}}Exponential\\ \end{tabular}} & \multirow{3}{*}{0.2259}                                                                                & \multirow{3}{*}{\textless{}0.0010}           & Max$\{Y_1,Y_2\}$          & 0.1803                                                                                       & 0.0030                                       \\
		&                                                                                    &                                                                                                        &                                              & Marginal of $Y_1$         & 0.1895                                                                                       & 0.0015                                       \\
		&                                                                                    &                                                                                                        &                                              & Marginal of $Y_2$         & 0.1551                                                                                       & 0.0162                                       \\ \hline
		\multirow{3}{*}{\begin{tabular}[c]{@{}l@{}}\textit{BPRHM1}\\ (R(1.2),1.3,1.2,1.0)\end{tabular}} & \multirow{3}{*}{\begin{tabular}[c]{@{}l@{}}Rayleigh\\ \end{tabular}}  & \multirow{3}{*}{0.2500}                                                                                & \multirow{3}{*}{\textless{}0.0010}           & Max$\{Y_1,Y_2\}$          & 0.1452                                                                                       & 0.0294                                       \\
		&                                                                                    &                                                                                                        &                                              & Marginal of $Y_1$         & 0.0878                                                                                       & 0.4242                                       \\
		&                                                                                    &                                                                                                        &                                              & Marginal of $Y_2$         & 0.1598                                                                                       & 0.0121                                       \\ \hline
		\multirow{3}{*}{\begin{tabular}[c]{@{}l@{}}\textit{BPRHM2}\\ (IW(1.2),1.2,1.4,1.6,1.8)\end{tabular}} & \multirow{3}{*}{\begin{tabular}[c]{@{}l@{}}Inverse\\ Weibull\end{tabular}}         & \multirow{3}{*}{0.1531}                                                                                      & \multirow{3}{*}{\textgreater 0.1000}                            & Max$\{Y_1,Y_2\}$          & 0.0553                                                                                       & 0.9196                                       \\
		&                                                                                    &                                                                                                        &                                              & Marginal of $Y_1$         & 0.0826                                                                                       & 0.5020                                       \\
		&                                                                                    &                                                                                                        &                                              & Marginal of $Y_2$         & 0.1177                                                                                       & 0.1254                                       \\ \hline
		\multirow{3}{*}{\begin{tabular}[c]{@{}l@{}}\textit{BPRHM2}\\ (E(2.5),1.2,1.4,1.6,1.8)\end{tabular}} & \multirow{3}{*}{\begin{tabular}[c]{@{}l@{}}Exponential\\ \end{tabular}} & \multirow{3}{*}{0.8660}                                                                  & \multirow{3}{*}{\textless 0.0010}       & Max$\{Y_1,Y_2\}$          & 0.7576                                                                                       & \textless 2.2 $\times 10^{-16}$                       \\
		&                                                                                    & \multicolumn{1}{l}{}                                                                                   & \multicolumn{1}{l}{}                         & Marginal of $Y_1$         & 0.7420                                                                                       & \textless 2.2 $\times 10^{-16}$                       \\
		&                                                                                    & \multicolumn{1}{l}{}                                                                                   & \multicolumn{1}{l}{}                         & Marginal of $Y_2$         & 0.6935                                                                                       & \textless 2.2 $\times 10^{-16}$                      \\ \hline
		\multirow{3}{*}{\begin{tabular}[c]{@{}l@{}}\textit{BPRHM2}\\ (IE(2.0),1.2,1.4,1.6,1.8)\end{tabular}} & \multirow{3}{*}{\begin{tabular}[c]{@{}l@{}}Inverse\\ Exponential\end{tabular}}     & \multirow{3}{*}{0.3613}                                                                  & \multirow{3}{*}{\textless 0.0010}    & Max$\{Y_1,Y_2\}$          & 0.3668                                                                                      & 4.106 $\times 10^{-12}$                       \\
		&                                                                                    & \multicolumn{1}{l}{}                                                                                   & \multicolumn{1}{l}{}                         & Marginal of $Y_1$         & 0.2779                                                                                       & 3.923 $\times 10^{-7}$                       \\
		&                                                                                    & \multicolumn{1}{l}{}                                                                                   & \multicolumn{1}{l}{}                         & Marginal of $Y_2$         & 0.4235                                                                                      & 5.551 $\times 10^{-16}$ \\ \hline                    
	\end{tabular}%
}
\end{table}

\begin{table}
	\centering
	\caption{Goodness of fit test for different baseline distributions (20\% censoring)\label{tabks1}}
	\resizebox{1\columnwidth}{!}{%
		\begin{tabular}{lllcc}
			\hline
			Model                                     & Baseline distribution            &             Variable                     & K-S test statistic & P value \\ \hline 
			\multirow{3}{*}{\textit{BPRHM2}(IW(2.1),1.5,1.6,2.0,1.8)} & \multirow{3}{*}{Inverse Weibull} & Max\{$Y_1,Y_2$\} & 0.0731             & 0.6596  \\
			&                                  & Marginal of $Y_1$                & 0.1292             & 0.0708  \\
			&                                  & Marginal of $Y_2$                & 0.1004             & 0.2661 \\ \hline 
			\multirow{3}{*}{\textit{BPRHM2}(E(2.0),1.5,1.6,2.0,1.8)} & \multirow{3}{*}{Exponential} & Max\{$Y_1,Y_2$\} & 0.6865             &  $< 2.2 \times 10^{-16}$  \\
			&                                  & Marginal of $Y_1$                & 0.7473             & $< 2.2 \times 10^{-16}$ \\
			&                                  & Marginal of $Y_2$                & 0.6652            & $< 2.2 \times 10^{-16}$ \\ \hline 
			\multirow{3}{*}{\textit{BPRHM2}(IE(1.8),1.5,1.6,2.0,1.8)} & \multirow{3}{*}{Inverse Exponential} & Max\{$Y_1,Y_2$\} & 0.6475             & $< 2.2 \times 10^{-16}$ \\
			&                                  & Marginal of $Y_1$                & 0.5353            & $< 2.2 \times 10^{-16}$  \\
			&                                  & Marginal of $Y_2$                & 0.5825             & $< 2.2 \times 10^{-16}$ \\ \hline
		\end{tabular}%
	}
\end{table}
\section{Application and Conclusion}\label{realdata}
In many situations, we may have to generate pairs of random variables from a continuous bivariate distribution. There exist specific simulation methods for certain bivariate distributions such as bivariate gamma, exponential and normal. However, in most scenarios we would run into multiple challenges while simulating from a bivariate distribution. The characterization result in Theorem \ref{ch1} helps in simulating a pair of random variables with distributions in the PRHR model. Thus, we can test for the goodness of fit of a bivariate data set with ease as it requires only the testing of univariate quantities.
\par We considered the  American Football dataset of the National Football League [\citet{csorgo1989testing}]. Here, $Y_1$ represents the the game time to the first points scored by kicking the ball between goal posts and $Y_2$ represents the game time by moving the ball into the end zone. There are 42 pairs of observations. The data points have been converted to decimal minutes similarly as in \citet{csorgo1989testing} and are divided by $100$ for computational purposes. \citet{kundu2010class} has analysed this dataset using \textit{BPRHM1} model with baseline distributions as Weibull, Rayleigh, Exponential and Linear Failure Rate. They estimated the parameters for the model with different baseline distributions. We used these estimates to calculate the univariate K-S test statistics and plot the theoretical distribution functions of each model. The Akaike Information Criterion (AIC) values are calculated for each model and is given in Table \ref{tabksdata}. 

\begin{table}[h!]
	\centering
	\caption{Goodness of fit test for American Football dataset under the \textit{BPRHM1} model \label{tabksdata}}
	\begin{tabular}{lcccccccc}
		\hline 
		\begin{tabular}[c]{@{}l@{}}Baseline \\ Distribution\end{tabular} & \multicolumn{2}{c}{\begin{tabular}[c]{@{}c@{}}Exponential \\ \end{tabular}} & \multicolumn{2}{c}{\begin{tabular}[c]{@{}c@{}}Weibull \\ \end{tabular}} & \multicolumn{2}{c}{\begin{tabular}[c]{@{}c@{}}Rayleigh \\ \end{tabular}} & \multicolumn{2}{c}{\begin{tabular}[c]{@{}c@{}}Linear \\ Failure Rate\end{tabular}}   \\ \hline
		AIC                                                              & \multicolumn{2}{c}{84.5}                                                               & \multicolumn{2}{c}{55.12}                                                            & \multicolumn{2}{c}{81.06}                                                             & \multicolumn{2}{c}{79.5}                                                             \\ \hline
		Variable                                                         & \begin{tabular}[c]{@{}c@{}}K-S test\\ statistic\end{tabular}   & P value               & \begin{tabular}[c]{@{}c@{}}K-S test \\ statistic\end{tabular} & P value              & \begin{tabular}[c]{@{}c@{}}K-S test \\ statistic\end{tabular}  & P value              & \begin{tabular}[c]{@{}c@{}}K-S test \\ statistic\end{tabular} & P value    \\ \hline
		Max\{$Y_1,Y_2$\}                                 & 0.1343                                                                  & 0.4000            & 0.1137                                                                & 0.6093         & 0.2111                                                                  & 0.0403          & 0.1644                                                                 & 0.1847       \\ \hline
		Marginal of $Y_1$                                                                 & 0.1771                                                                  & 0.1435            & 0.1584                                                                 & 0.2426         & 0.2387                                                                  & 0.0167         & 0.2035                                                                & 0.0616         \\ \hline
		Marginal of $Y_2$                                                                 & 0.1503                                                                  & 0.2709          & 0.1300                                                                  & 0.4402          & 0.2262                                                                 & 0.0228          & 0.1800                                                                & 0.1160          \\ \hline
	\end{tabular}
\end{table}

\begin{figure}
	\begin{minipage}{0.325\textwidth}
		\includegraphics[width=1\textwidth]{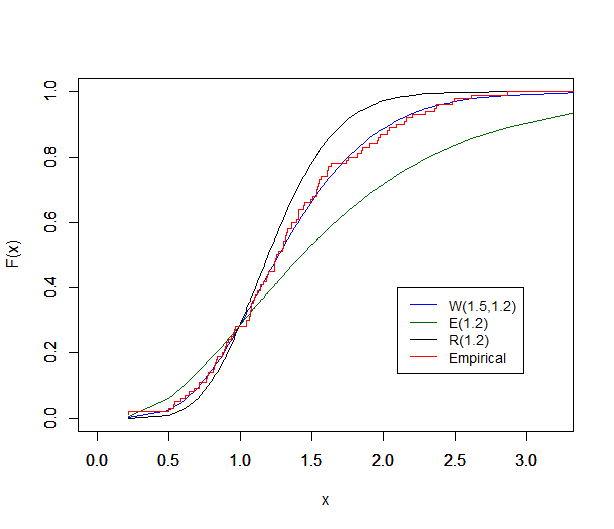} 
		\subcaption{Max\{$Y_1,Y_2$\}}
	\end{minipage}\hfill
	\begin{minipage}{0.325\textwidth}
		\includegraphics[width=1\textwidth]{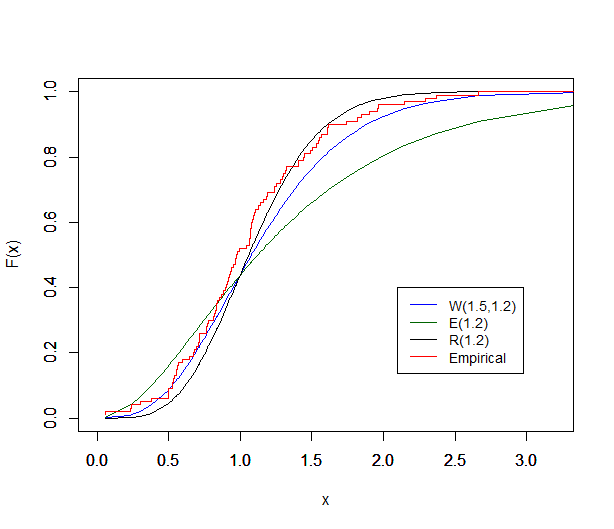} 
		\subcaption{Marginal of $Y_1$}
	\end{minipage}
	\begin{minipage}{0.325\textwidth}
		\includegraphics[width=1\textwidth]{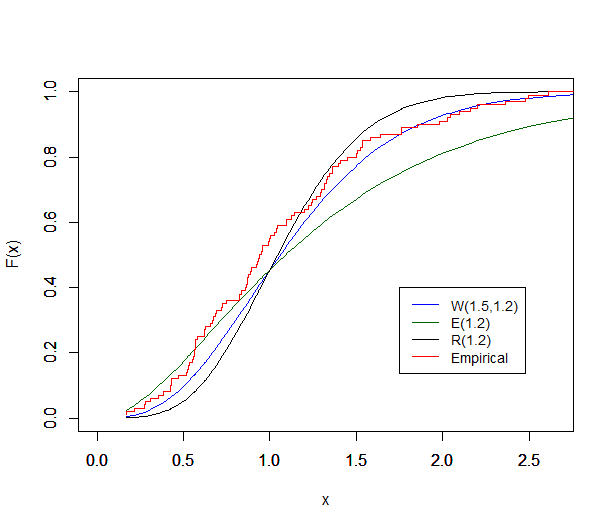} 
		\subcaption{Marginal of $Y_2$}
	\end{minipage}	
	\caption{Plots of theoretical and empirical distribution functions of \textit{BPRHM1} model ($BPRHM1(F_{0},1.3,1.2,1.0)$)}\label{figBPRHsim}
	\begin{minipage}{0.325\textwidth}
		\includegraphics[width=1\textwidth]{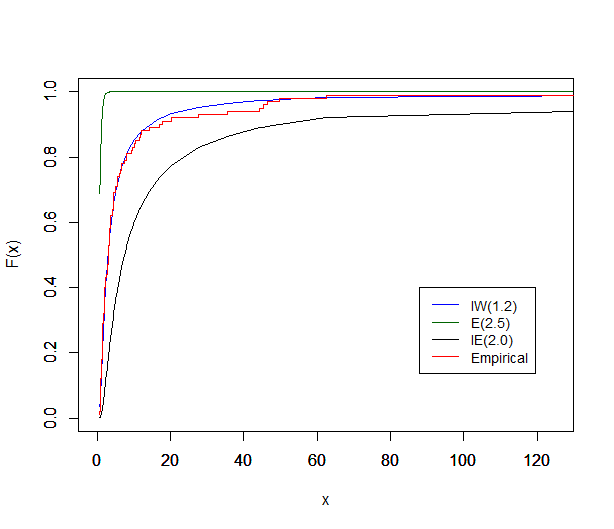} 
		\subcaption{Max\{$Y_1,Y_2$\}}
	\end{minipage}\hfill
	\begin{minipage}{0.325\textwidth}
		\includegraphics[width=1\textwidth]{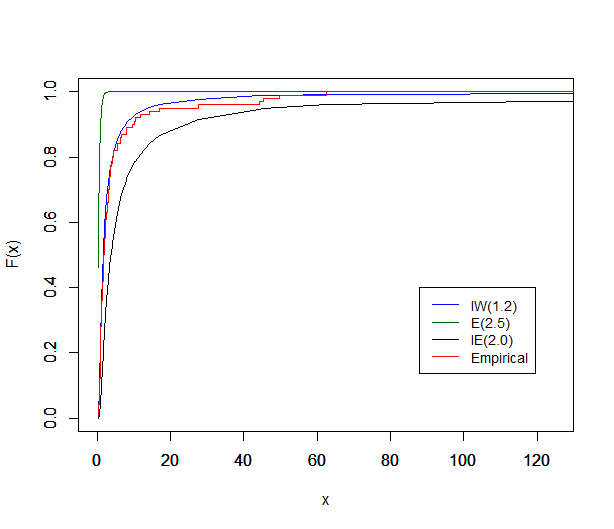} 
		\subcaption{Marginal of $Y_1$}
	\end{minipage}
	\begin{minipage}{0.325\textwidth}
		\includegraphics[width=1\textwidth]{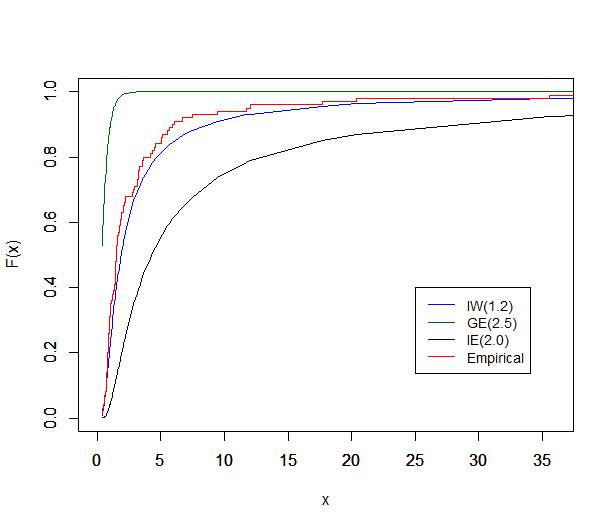} 
		\subcaption{Marginal of $Y_2$}
	\end{minipage}	
	\caption{Plots of theoretical and empirical distribution functions of \textit{BPRHM2} model ($BPRHM2(F_{0},1.2,1.4,1.6,1.8)$)}\label{figDPRHsim}
	\begin{minipage}{0.325\textwidth}
		\includegraphics[width=1\textwidth]{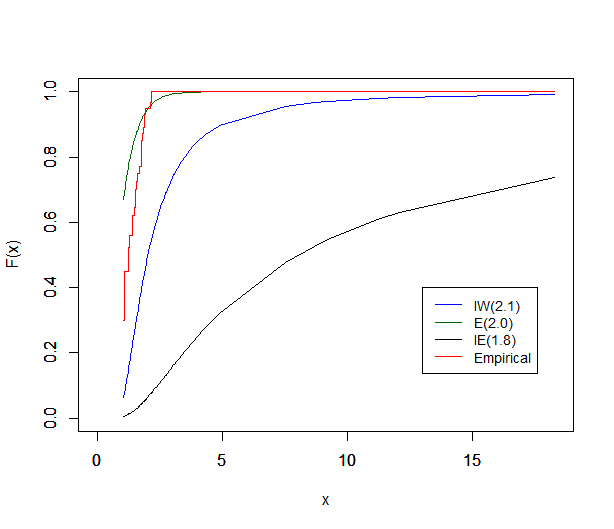} 
		\subcaption{Max\{$Y_1,Y_2$\}}
	\end{minipage}\hfill
	\begin{minipage}{0.325\textwidth}
		\includegraphics[width=1\textwidth]{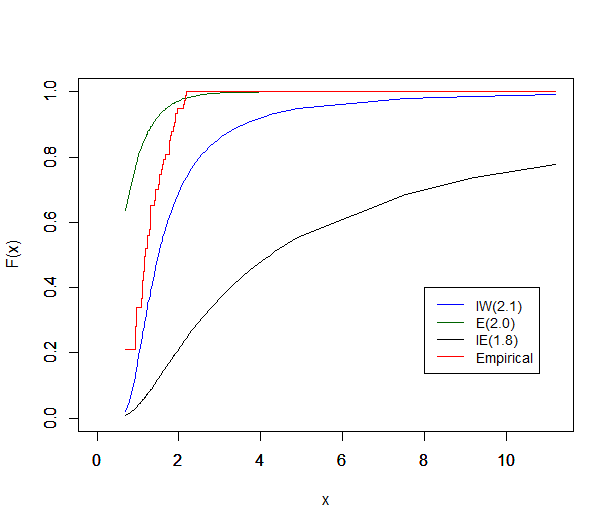} 
		\subcaption{Marginal of $Y_1$}
	\end{minipage}
	\begin{minipage}{0.325\textwidth}
		\includegraphics[width=1\textwidth]{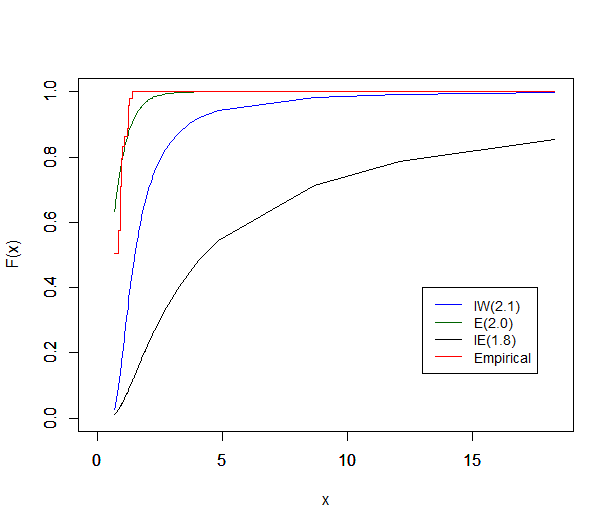} 
		\subcaption{Marginal of $Y_2$}
	\end{minipage}	
	\caption{Plots of theoretical and empirical distribution functions of \textit{BPRHM2} model with 20\% censoring ($BPRHM2(F_{0},1.5,1.6,2.0,1.8)$)}\label{figDPRHsimcens}
	\begin{minipage}{0.325\textwidth}
		\includegraphics[width=1\textwidth]{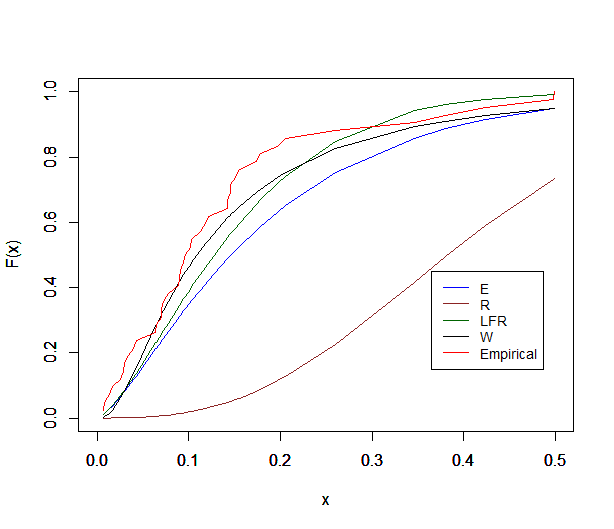} 
		\subcaption{Max\{$Y_1,Y_2$\}}
	\end{minipage}\hfill
	\begin{minipage}{0.325\textwidth}
		\includegraphics[width=1\textwidth]{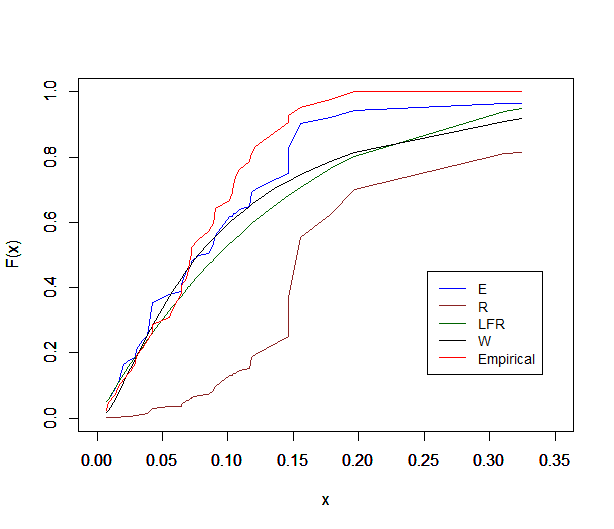} 
		\subcaption{Marginal of $Y_1$}
	\end{minipage}
	\begin{minipage}{0.325\textwidth}
		\includegraphics[width=1\textwidth]{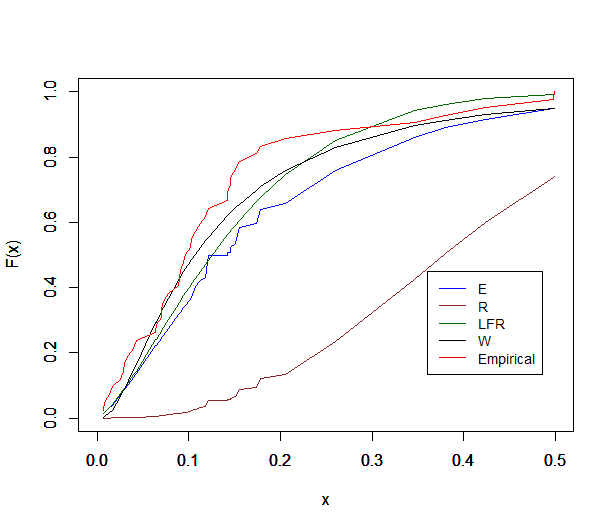} 
		\subcaption{Marginal of $Y_2$}
	\end{minipage}	
	\caption{Plots of theoretical and empirical distribution functions of \textit{BPRHM1} model with different baseline distributions for the American Football dataset}\label{figBPRHdata}
\end{figure}

It is noted that the AIC value is lowest $(55.12)$ for the model with Weibull as baseline distribution which makes it a suitable model for the American Football dataset. As seen from Table \ref{tabksdata}, the p-values of the univariate K-S tests are also highest for this model. Hence, when the AIC values are used for model selection, we can also append the goodness of fit test results to it. Figure \ref{figBPRHdata} shows the theoretical and empirical plots of the different variables under this model with different baseline distributions. Thus, the findings are consistent with the results of \citet{kundu2010class}.
\subsection*{\textbf{Funding}}
The authors disclosed the receipt of the following financial support for the research, authorship, and/or publication of this article. This work was supported by the Council of Scientific \& Industrial Research (CSIR) through the JRF scheme vide No 09/239(0551)/2019-EMR-I.

	
%
	\bibliography{ref}

\begin{thebibliography}{}

\bibitem[Arias-Nicol{\'a}s et~al., 2016]{arias2016new}
Arias-Nicol{\'a}s, J.~P., Ruggeri, F., and Su{\'a}rez-Llorens, A. (2016).
\newblock New classes of priors based on stochastic orders and distortion
  functions.
\newblock {\em Bayesian Analysis}, 11(4):1107--1136.

\bibitem[Asha and John, 2007]{asha2007models}
Asha, G. and John, R.~C. (2007).
\newblock Models characterized by the reversed lack of memory property.
\newblock {\em Calcutta Statistical Association Bulletin}, 59(1-2):1--14.

\bibitem[Barlow et~al., 1963]{barlow1963properties}
Barlow, R.~E., Marshall, A.~W., Proschan, F., et~al. (1963).
\newblock Properties of probability distributions with monotone hazard rate.
\newblock {\em The Annals of Mathematical Statistics}, 34(2):375--389.

\bibitem[Bismi and Ramachandran~Nair, 2005]{bismi2005bivarite}
Bismi, N.~G. and Ramachandran~Nair, V. (2005).
\newblock {\em Bivarite burr distributions}.
\newblock PhD thesis, Cochin University of Science and Technology.

\bibitem[Block et~al., 1998]{block1998reversed}
Block, H.~W., Savits, T.~H., and Singh, H. (1998).
\newblock The reversed hazard rate function.
\newblock {\em Probability in the Engineering and informational Sciences},
  12(1):69--90.

\bibitem[Cs{\"o}rg{\H o} and Welsh, 1989]{csorgo1989testing}
Cs{\"o}rg{\H o}, S. and Welsh, A. (1989).
\newblock Testing for exponential and marshall–olkin distributions.
\newblock {\em Journal of Statistical Planning and Inference}, 23(3):287--300.

\bibitem[Di~Crescenzo, 2000]{di2000some}
Di~Crescenzo, A. (2000).
\newblock Some results on the proportional reversed hazards model.
\newblock {\em Statistics \& probability letters}, 50(4):313--321.

\bibitem[Ganser and Hewett, 2010]{ganser2010accurate}
Ganser, G.~H. and Hewett, P. (2010).
\newblock An accurate substitution method for analyzing censored data.
\newblock {\em Journal of occupational and environmental hygiene},
  7(4):233--244.

\bibitem[Gupta et~al., 1998]{gupta1998modeling}
Gupta, R.~C., Gupta, P.~L., and Gupta, R.~D. (1998).
\newblock Modeling failure time data by lehman alternatives.
\newblock {\em Communications in Statistics-Theory and methods},
  27(4):887--904.

\bibitem[Gupta and Gupta, 2007]{gupta2007proportional}
Gupta, R.~C. and Gupta, R.~D. (2007).
\newblock Proportional reversed hazard rate model and its applications.
\newblock {\em Journal of Statistical Planning and Inference},
  137(11):3525--3536.

\bibitem[G{\"u}rler, 1996]{gurler1996bivariate}
G{\"u}rler, {\"U}. (1996).
\newblock Bivariate estimation with right-truncated data.
\newblock {\em Journal of the American Statistical Association},
  91(435):1152--1165.

\bibitem[Justel et~al., 1997]{justel1997multivariate}
Justel, A., Pe{\~n}a, D., and Zamar, R. (1997).
\newblock A multivariate kolmogorov-smirnov test of goodness of fit.
\newblock {\em Statistics \& Probability Letters}, 35(3):251--259.

\bibitem[Keilson and Sumita, 1982]{keilson1982uniform}
Keilson, J. and Sumita, U. (1982).
\newblock Uniform stochastic ordering and related inequalities.
\newblock {\em Canadian Journal of Statistics}, 10(3):181--198.

\bibitem[K{\i}z{\i}laslan, 2017]{kizilaslan2017bayesian}
K{\i}z{\i}laslan, F. (2017).
\newblock The e-bayesian and hierarchical bayesian estimations for the
  proportional reversed hazard rate model based on record values.
\newblock {\em Journal of Statistical Computation and Simulation},
  87(11):2253--2273.

\bibitem[Krishnamoorthy et~al., 2009]{krishnamoorthy2009model}
Krishnamoorthy, K., Mallick, A., and Mathew, T. (2009).
\newblock Model-based imputation approach for data analysis in the presence of
  non-detects.
\newblock {\em Annals of Occupational Hygiene}, 53(3):249--263.

\bibitem[Kundu et~al., 2014]{kundu2014multivariate}
Kundu, D., Franco, M., and Vivo, J.-M. (2014).
\newblock Multivariate distributions with proportional reversed hazard
  marginals.
\newblock {\em Computational Statistics \& Data Analysis}, 77:98--112.

\bibitem[Kundu and Gupta, 2004]{kundu2004characterizations}
Kundu, D. and Gupta, R.~D. (2004).
\newblock Characterizations of the proportional (reversed) hazard model.
\newblock {\em Communications in Statistics - Theory and Methods},
  33(12):3095--3102.

\bibitem[Kundu and Gupta, 2010]{kundu2010class}
Kundu, D. and Gupta, R.~D. (2010).
\newblock A class of bivariate models with proportional reversed hazard
  marginals.
\newblock {\em Sankhya B}, 72(2):236--253.

\bibitem[Lawless, 2011]{lawless2011statistical}
Lawless, J.~F. (2011).
\newblock {\em Statistical models and methods for lifetime data}, volume 362.
\newblock John Wiley \& Sons.

\bibitem[Mudholkar and Srivastava, 1993]{mudholkar1993exponentiated}
Mudholkar, G.~S. and Srivastava, D.~K. (1993).
\newblock Exponentiated weibull family for analyzing bathtub failure-rate data.
\newblock {\em IEEE transactions on reliability}, 42(2):299--302.

\bibitem[Mudholkar et~al., 1995]{mudholkar1995exponentiated}
Mudholkar, G.~S., Srivastava, D.~K., and Freimer, M. (1995).
\newblock The exponentiated weibull family: A reanalysis of the
  bus-motor-failure data.
\newblock {\em Technometrics}, 37(4):436--445.

\bibitem[Popovi{\'c} et~al., 2021]{popovic2021generalized}
Popovi{\'c}, B.~V., Gen{\c{c}}, A.~{\.I}., and Domma, F. (2021).
\newblock Generalized proportional reversed hazard rate distributions with
  application in medicine.
\newblock {\em Statistical Methods \& Applications}, pages 1--22.

\bibitem[Roy, 2002]{roy2002characterization}
Roy, D. (2002).
\newblock A characterization of model approach for generating bivariate life
  distributions using reversed hazard rates.
\newblock {\em Journal of the Japan Statistical Society}, 32(2):239--245.

\bibitem[Ruggeri et~al., 2021]{ruggeri2021new}
Ruggeri, F., S{\'a}nchez-S{\'a}nchez, M., Sordo, M.~{\'A}., and
  Su{\'a}rez-Llorens, A. (2021).
\newblock On a new class of multivariate prior distributions: Theory and
  application in reliability.
\newblock {\em Bayesian Analysis}, 16(1):31--60.

\bibitem[Sankaran and Gleeja, 2008]{sankaran2008proportional}
Sankaran, P. and Gleeja, V. (2008).
\newblock Proportional reversed hazard and frailty models.
\newblock {\em Metrika}, 68(3):333--342.

\bibitem[Vasudevan and Asha, 2021]{durga2021proportional}
Vasudevan, D. and Asha, G. (2021).
\newblock A proportional reversed hazards model for load sharing systems.
\newblock {\em Personal communication}.

\bibitem[Ware and Demets, 1976]{ware1976reanalysis}
Ware, J.~H. and Demets, D.~L. (1976).
\newblock Reanalysis of some baboon descent data.
\newblock {\em Biometrics}, pages 459--463.

\end{thebibliography}

\end{document}